\documentclass[12pt]{article}

\usepackage{amsmath,amsfonts,amssymb,amsthm,amscd,enumerate}
\usepackage[all]{xy}
\usepackage{graphicx}
\usepackage{epstopdf}

\usepackage[dvipsnames]{color}

\title{Closed quantum surfaces from the Toeplitz extension}

\author{{\sc Arley Sierra} \\  
\footnotesize
Centro de Ciencias Matem\'aticas, Campus Morelia\\
\footnotesize
Universidad Nacional Aut\'onoma de M\'exico (UNAM), Morelia, M\'exico\\
\footnotesize
E-mail: {\it arleysierra23@gmail.com}\\[12pt] 
{\sc Elmar Wagner} \\
\footnotesize
Instituto de F\'isica y Matem\'aticas\\
\footnotesize
Universidad Michoacana de San Nicol\'as de Hidalgo, Morelia, M\'exico\\
\footnotesize
E-mail: {\it elmar@ifm.umich.mx}}

\date{}                                           

\newtheorem{thm}{Theorem}

\theoremstyle{definition}
\newtheorem{defn}[thm]{Definition}

\newcommand{\nc}[2]{\newcommand{#1}{#2}}
\newcommand{\rnc}[2]{\renewcommand{#1}{#2}}
\rnc{\[}{\begin{equation}}
\rnc{\]}{\end{equation}}
\nc{\wegengruen}{\end{equation}}

\newcommand{\Z}{\mathbb{Z}}
\newcommand{\N}{\mathbb{N}}
\newcommand{\R}{\mathbb{R}}
\newcommand{\C}{\mathbb{C}}

\newcommand{\bT}{\mathbb{T}}
\newcommand{\bS}{\mathbb{S}^1}

 \newcommand{\D}{\mathbb{D}}  
\newcommand{\hsp}{{\hspace{-1pt}}}
\newcommand{\hs}{{\hspace{1pt}}}

\newcommand{\cA}{{\mathcal{A}}}
\newcommand{\cK}{{\mathcal{K}}} 
\newcommand{\cB}{{\mathcal{B}}}
\newcommand{\cM}{{\mathcal{M}}}
 \newcommand{\T}{\mathcal{T}}  
 \newcommand{\cC}{\mathcal{C}}
  
  \newcommand{\Cal}{\mathfrak{C}}

\newcommand{\lN}{\ell_2(\N_0)}
\newcommand{\lZ}{\ell_2(\Z)}
\def\CS{\cC(\mathbb{S}^1)} 

\def\bD{\bar\D}
\def\LD{L_2(\D)}
\def\AD{A_2(\D)} 
\def\KlN{\cK(\lN)}

\def\CD{\cC(\bar\D)}
\def\CTq{\cC(\mathbb{T}^g_{q})} 
\def\CPq{\cC(\mathbb{P}^n_{k,q})}
\def\CSq{\cC(\mathbb{S}_{q}^2)} 

\def\s{\sigma} 
\def\a{\alpha}

\newcommand{\im}{\mathrm{i}}
\newcommand{\e}{\mathrm{e}}

\newcommand{\id}{{\mathrm{Id}}}
\newcommand{\ind}{{\mathrm{ind}}}

\newcommand{\spec}{\mathrm{spec}}

\newcommand{\wind}{\mathrm{wind}}

\newcommand{\ra}{\rightarrow}
\newcommand{\lra}{\longrightarrow}

\newcommand{\mplus}[2]{\underset{{#1}}{\overset{{#2}}{\mbox{$\oplus$}}}}
\newcommand{\mfrac}[2]{\mbox{$\frac{#1}{#2}$}}

\begin{document}
\maketitle

\begin{abstract} 
Closed quantum surfaces of any genus are defined as subalgebras of the Toeplitz algebra 
by mimicking the classical construction of identifying arcs on the boundary 
of the (quantum) unit disk. Isomorphism classes obtained from different arrangements of arcs are classified. 
It is shown that the K-groups are isomorphic to the classical counterparts and explicit generators 
of the C*-algebras and of the K-groups are given. 
\end{abstract}

\section{Introduction} 
Quantization of a compact topological space or manifold means, roughly speaking, 
the replacement of the C*-algebra of continuous functions by a 
noncommutative C*-algebra. However, there is no universal procedure that tells us how to pass from a commutative C*-algebra 
to a noncommutative one while maintaining certain topological features of the space. This reflects a recurrent problem in quantum physics, 
where no functorial method for the quantization of classical observables or fields is known \cite{G}. On the other hand, 
the relevance of noncommutative geometry \cite{C} in mathematics and theoretical physics can only be evidenced by providing 
a proper amount of useful examples. 
The aim of this paper is to present a whole family of noncommutative topological spaces, namely quantizations of all closed 
two-dimensional surfaces. This can been seen as a first step of the wider project of quantizing (finite) CW-com\-plexes~\cite{qCW}. 

Our staring point will be the Toeplitz quantization of the unit disk \cite{KL}. It replaces the continuous functions on the closed disk 
by their corresponding Toeplitz operators, see Section \ref{STA} for more details. 
The C*-algebra $\T$ generated by all these Toeplitz operators yields a non-trivial 
C*-algebra extension of $\CS$ by the compact operators $\cK$ 
on a separable Hilbert space. Then the so-called symbol map $\s: \T\ra \T/\cK \cong \CS$ may be viewed as the restriction 
of continuous functions on the quantum disk to the boundary circle $\bS$.  Based on the premise that the quantum disk admits 
a classical boundary circle, we will define in Section \ref{Sdef} closed quantum surfaces by considering C*-subalgebras of $\T$ 
that correspond to glueing pairs of arcs on the boundary circle in such a way that the commutative analog yields 
a C*-algebra isomorphic to the continuous functions on a specific closed surface. 
 
Classically, a closed surface can be obtained from different arrangements of arcs. The standard method for providing an 
homeomorphism is based on a ``cut and glue'' technique, which is not available in the quantum case. In fact, we will show in 
Section \ref{Siso} that isomorphism classes of closed quantum surfaces are labeled by the number of projective spaces and tori 
that are used when the arrangements of arcs are directly interpreted as a connected sum of these building blocks. In this sense, 
the Toeplitz quantization decreases degeneracy. Furthermore, for each isomorphism class, 
we will give a essentially normal generator of the corresponding C*-algebra in terms of unilateral and bilateral shift operators. 

Since the definition of closed quantum surfaces will be given just by analogy to the classical case, there arises the question 
of whether the quantization changes topological invariants. The topological invariants that we consider in this paper are the 
K-groups of the C*-algebras. In Section \ref{SKT}, we will prove that the K-groups of the closed quantum surfaces 
are isomorphic to the classical counterparts. Finally, for eventual future use, explicit descriptions of the generators 
of the K-groups are given.

\section{Quantum disk view on the Toeplitz algebra}    \label{STA} 

Let $\D:= \{ z\in \C: |z|<1\}$ and $\bD:= \{ z\in \C: |z|\leq 1\}$ denote the open and closed unit disk, respectively. 
We write $\LD$ for  the Hilbert space of square-integrable functions with respect to 
the Lebesgue measure and  $\AD$ for the subspace 
of square-integrable holomorphic functions on $\D$. Since $\AD\subset \LD$ is closed, there 
exists an orthogonal projection, say $P$, from $\LD$ onto $\AD$. 
Now the Toeplitz operator $T_f\in \cB(\AD)$ with continuous symbol $f\in \CD$ is given by 
$$
 T_f(\psi):= P(f\hs \psi), \quad \psi \in \AD \subset \LD, 
$$
and the Toeplitz algebra $\T$ may be defined as the C*-subalgebra generated by all $T_f$ in the 
C*-algebra of bounded operators $\cB:= \cB(\AD)$.  

It can be shown (see e.\,g.\ \cite{V}) that the operador ideal of compact operators $\cK:= \cK(\AD) \cong \KlN$ 
belongs to $\T$ and that the quotient $\T/\cK$ is isomorphic to $\CS $, where we view $\mathbb{S}^1= \partial\bD$ as 
the boundary of $\bD$. This gives rise to the C*-algebra extension 
\[ \label{Cex1} 
\xymatrix{
 0\;\ar[r]&\; \cK\;\ar[r]^{\iota\ \ } &\; \T \;\ar[r]^{\s\ \ } &\; \CS\;\ar[r] &\ 0\, ,} 
 \]
with the so-called symbol map $\s : \T \lra \CS$ given by $\s(T_f) =  f\!\! \upharpoonright_{\mathbb{S}^1}$ 
for all $f\in \CD$. 

The application $\CD \ni f \mapsto T_f\in\cB(\AD)$  will be viewed as a quantization  
of the commutative unital C*-algebra $\cC(\bD)$. In agreement with \cite{KL}, we refer 
to the Toeplitz algebra $\T=:\cC(\bD_q)$ as the 
algebra of continuous functions on the quantum disk. 
In the commutative case, the C*-algebra extension \eqref{Cex1} corresponds to the exact sequence 
\[ \label{Cex2} 
\xymatrix{
 0\;\ar[r]&\; \cC_0(\D)\;\ar[r]^{\iota\ \ } &\; C(\bD) \;\ar[r]^{\rho\ \ } &\; \CS\;\ar[r] &\ 0\, ,} 
 \]
where $\rho(f)= f\!\! \upharpoonright_{\mathbb{S}^1}$. 

Let $z\in\CD$, $z(x)=x$ denote the identity function. By the Stone--Weierstrass Theorem, 
the functions $1$, $z$ and $z^*:=\bar z$ generate the C*-algebra $\CD$. For this reason, 
$1$, $T_z$ and $T_{\bar z}$ generate $\T$. On the orthonormal basis 
$\{ e_n:=\frac{\sqrt{n+1}}{\sqrt{\pi}} \hs z^n : n\in\N_0\}$ of $\AD$, the operator $T_z$ acts by 
$T_z e_n = \frac{\sqrt{n+1}}{\sqrt{n+2}}e_{n+1}$. Next, consider the  
unilateral shift 
\[ \label{S}
Se_n := e_{n+1}, \quad n\in \N_0. 
\] 
As $\underset{n\ra\infty}{\lim} \frac{\sqrt{n+1}}{\sqrt{n+2}} -1 =0$, it follows that $T_z -S \in \cK$. 
Knowing that the C*-al\-ge\-bra generated by $1$, $S$ and $S^*$ contains the compact operators, 
it can be inferred that $1$, $S$ and $S^*$ also generate $\T$. Moreover, 
\[ \label{sS}
\s(S)= \s(T_z) =: u\in\CS, \quad u(\e^{\im t}) = \e^{\im t}. 
\] 

Comparing the C*-algebra extensions \eqref{Cex1} and \eqref{Cex2}, it seems that the 
Toeplitz quantization $f\mapsto T_f$ amounts to replacing $\cC_0(\D)$ by $\cK$. 
The following chain of K-theoretic identities gives an additional motivation 
for this interpretation: 
$$
K_i(\cC_0(\D)) \cong K_i(\Sigma^2\C) \cong K_i(\C) \cong K_i(\cK\otimes \C) \cong K_i(\cK), \quad i=0,1. 
$$
Here, $\Sigma\cA$ denotes the suspension of a C*-algebra $\cA$. The first isomorphism comes from 
an isomorphism of C*-algebras, the second one from Bott periodicity, the third one holds by stabilization, and the 
last one is trivial.  Using $K_0(\C)= \Z[1]$, $K_1(\C)= 0$, $K_0(\CS)= \Z[1]$ and $K_1(\CS)= \Z[u]$, with 
$u$ being the unitary defined in \eqref{sS}, the K-groups of $\T$ can easily be computed from 
the 6-term exact sequence of K-theory: 
\begin{equation}  \label{stX}
\xymatrixcolsep{2pc}
\xymatrix{  
 \Z[1\!-\!SS^*] \cong K_0(\cK) \ \ar[r]^-{\iota_\ast} & 
  \ K_0 ( \T)\   \ar[r]^-{{\s}_\ast} &
  \ K_0 (\CS)\cong \Z[1] \ \ar[d]^{\mathrm{exp}}\\
\Z[u]\cong K_1 (\CS)\ \ar[u]^{\ind} &
  \ K_1 ( \T)\  \ar[l]_-{\ \ {\s}_\ast} &
\ K_1 (\cK) \ar[l]_-{\ \ \iota_\ast}  \cong 0 \,. 
  }
\end{equation} 
The index map $\ind : K_1 (\CS) \ra \Z \cong K_0(\cK)$ relates closely to the 
winding number $\wind(\Phi)\in\Z$ of continuous functions  $\Phi:\ \mathbb{S}^1 \ra \mathbb{S}^1$ and to the 
Fredholm index $\mathrm{Ind}(F)\in\Z$ of invertible elements $[F]\in\Cal$ in the Calkin algebra $\Cal:=\cB/\cK$. 
Applying  the fact that $K_0(\cB) = K_1(\cB) = 0$, the index map $\ind : K_1(\Cal) \ra \Z \cong K_0(\cK)$ 
in the  6-term exact sequence corresponding to the extension $0\ra\cK\ra\cB\ra\Cal\ra 0$ 
yields an isomorphism. Moreover, the unitary element $u = \s(S)\in \Cal$ 
admits obviously a lift by the isometry $S$ from \eqref{S}. By \cite[Ex.\ 8.C]{WO}, 
$\ind[u] = -[1\!-\!SS^*]$, where $[1\!-\!SS^*]$ represents a generator of 
$K_0(\cK)$. Since  $\ind$ is a group homomorphism, we may write, for all $k\in\Z$, 
\[ \label{Sind} 
k=\wind[u^k] =-\ind[u^k] 
=\left\{ \begin{array}{ll} -\hs \mathrm{Ind}(S^{*|k|}),  & k<0,\\ -\hs \mathrm{Ind}(S^{k}),  & k\geq0,\end{array}\right.
\]
and since the K-theory of C*-algebras 
is homotopy invariant, we get 
\[ \label{wind} 
\ind[\Phi]  = \mathrm{Ind}(F_\Phi)  = - \wind[\Phi] 
\]
for any invertible function $\Phi\in\CS\cong \T/\cK \subset \Cal$, where $F_\Phi$ denotes a lift of~$\Phi$.  

There is also an analogy to the famous Bott projection of $\cC_0(\D)$, which illustrates nicely 
the interpretation of $\T$ as a quantization of $\CD$. 
Given an unitary function $\upsilon\in\CS$, 
let $\zeta := r\upsilon\in\CD$ be an extension to an continuous function on the closed disk, 
where $r$ denotes the radius function of the points in $\bD$. Then, by \cite[Ex.\ 8.D]{WO}, 
$$
\ind[\upsilon]= \left[\begin{pmatrix} \zeta\bar\zeta  & \zeta \sqrt{1-\bar\zeta\zeta} \\[2pt] 
                                                          \sqrt{1-\bar\zeta\zeta}\,\bar\zeta   & 1- \bar\zeta\zeta\end{pmatrix} \right] 
 - \left[\begin{pmatrix} 1 &  0 \\ 0 & 0\end{pmatrix} \right] \in K_0( \cC_0(\D))                                                        
$$
and 
\[ \label{PTz} 
\ind[\upsilon]= \left[\begin{pmatrix} T_\zeta T_\zeta^*  & T_\zeta \sqrt{1-T_\zeta^*T_\zeta} \\[2pt]
                                                          \sqrt{1-T_\zeta^*T_\zeta}\,T_\zeta^*   & 1- T_\zeta^*T_\zeta\end{pmatrix} \right] 
 - \left[\begin{pmatrix} 1 &  0 \\ 0 & 0\end{pmatrix} \right] \in K_0( \cK).                                                        
\]
Note the striking similarity between these projections. For $\upsilon = u\in\CS$ and $\zeta =ru=z\in\CD$ (the identity function on $\bD$), 
the first formula renders the Bott projection of $\cC_0(\D)$.

\section{Definition of closed quantum surfaces} \label{Sdef}

Classical closed surfaces (compact and without boundary) can be described by simply connected polygons in the 2-dimensional plane 
with a prescribed identification of the boundary edges \cite{F}. There is no loss in generality if we replace the 
polygon by the closed unit disk $\bD\subset \C$ and turn the edges into arcs on the boundary circle maintaining 
 their orientations. These arcs may be labeled by pairs of letters $a_1, a_2, \ldots$ and  $b_1, b_2, \ldots$ if they have the same 
orientation, or by pairs of letters $a_1, a_2, \ldots$ and $a_{k+1}^{-1}, a_{k+2}^{-1}, \ldots$ if they are given the opposite orientation. 
Assume that these arcs are parametrized by continuous curves on the interval $[0,1]$, e.g.\ $[0,1] \ni t \mapsto a_j(t) \in \partial\bD$, always 
in the direction of their orientation. Then glueing pairs of the arcs means identifying the points $a_j(t)$ and $b_j(t)$ 
 if two numbered sets of letters correspond to each other, or the points $a_j(t)$ and $a_j^{-1}(t)$ if the arc $a_j$ 
 occurs exactly once with its negative orientation $a_j^{-1}$ and has thus no companion  $b_j$.
 
 In this paper, only the following arrangements will be considered. 
 Given $g\in \N$, we use the notation $\bT^g$ if the boundary 
 $\partial \bD\cong \bS$ is divided into $4g$ arcs $a_1,\ldots, a_{2g}, a_1^{-1}, \ldots, a_{2g}^{-1}$ 
 such that the topological quotient $ \bar\D\,/\!\hsp\sim$ under the equivalence relations 
 \[ \label{Tsim} 
 z\sim z,\ \,\forall z\in\bD, \ \ a_j(t)\sim a_j^{-1}(t),  \  \,j=1,\ldots,2g, \ \,t\in[0,1], 
\]
yields a realization of a closed orientable surface of genus $g$. 

For $k,n\in \N$ with $k\leq n$, we write $\mathbb{P}^n_k$ for a division of the boundary $\partial \bD\cong \bS$ 
into $2n$ arcs $a_1,\ldots, a_{k}, b_1, \ldots, b_{k}, a_{k+1}, \ldots, a_{n},a_{k+1}^{-1}, \ldots, a_{n}^{-1}$ 
such that the topological quotient $\bar\D\,/\!\hsp\sim$\,,  
 \[ \label{Psim} 
  z\sim z,\ \,\forall z\in\bD, \ \   a_i(t)\sim b_i(t),\ \, a_j(t)\sim a_j^{-1}(t), \  \,i\leq k, \ \, j>k,  \ \,t\in[0,1], 
\]
is homeomorphic to a closed non-orientable surface of Euler genus $n$. 
Shrinking the arcs $a_{k+1}, \ldots a_{n}$, $a_{k+1}^{-1}, \ldots a_{n}^{-1}$ to a point 
yields a closed non-orientable surface of Euler genus $n-k$, 
whereas shrinking the arcs $a_1,\ldots, a_{k}$, $b_1, \ldots, b_{k}$ to a point yields an 
orientable surface of genus $(n-k)/2$. Obviously, $n-k$ has to be an even number. In this case, 
the equivalence relation in \eqref{Psim} corresponds to the connected sum 
\[ \label{TP}
 \mathbb{P}^n_k:=\bar\D\,/\!\sim\, \ \cong\  \underbrace{  \mathbb{P}^1_1 \#\cdots  \#\mathbb{P}^1_1}_{\substack{k\text{ times}}} 
 \,\#\, \underbrace{  \mathbb{T}^1 \#\cdots  \#\mathbb{T}^1}_{\substack{(n-k)/2\text{ times}}} 
 \ \cong\   \mathbb{P}^k_k \,\#\, \mathbb{T}^{(n-k)/2} , 
\] 
 which is known to be homeomorphic to the closed non-orientable surface $\mathbb{P}^n := \mathbb{P}^n_n$ 
 of Euler genus $n$. 
 
 The C*-algebras of continuous functions $\cC(\bT^g)$ and $\cC(\mathbb{P}^n_k)$ 
 on the surfaces  $\bT^g$ and $\mathbb{P}^n_k$ are then isomorphic to 
 the respective subalgebras of all functions $f\in C(\bD)$ such that $f(x)=f(y)$ whenever $x\sim y$, where the equivalence relations are 
 given in \eqref{Tsim} and \eqref{Psim}, respectively. Comparing the C*-algebra extensions \eqref{Cex1} and \eqref{Cex2}, 
 and viewing the symbol map $\s:\T\ra \CS$ as the restriction of continuous functions on the quantum disk to the boundary, 
 the next definition of closed quantum surfaces becomes fairly obvious. 
 
 \begin{defn} \label{D1} 
 For $g\in \N$, let the boundary $\partial \bD\cong \bS$ be divided into $4g$ arcs $a_1,\ldots, a_{2g}, a_1^{-1}, \ldots, a_{2g}^{-1}$ 
 such that the topological quotient $\bT^g := \bar\D\,/\!\hsp\sim$ with the equivalence relation given in \eqref{Tsim} 
yields a realization of a closed orientable surface of genus $g$. Then an  
orientable closed quantum surface of genus $g$ 
is defined by the C*-algebra 
\[ \label{Tq} 
\CTq := \{ f\hsp\in\hsp\T: \s(f)(x)\hsp= \hsp\s(f)(y),  \  \,\forall x,y\in\partial \bD \,\text{ such that } x\sim y \}, 
\]
where $\s:\T\ra \CS$ denotes the symbol map. Likewise, a quantum 2-\-sphere is given by 
\[ \label{Sq} 
\CSq := \{ f \in \T: \s(f)( \e^{\pi\im t}) = \s(f)( \e^{-\pi\im t}), \ \,t\in[0,1]\}. 
\] 
Next, for $n\in \N$ and $k\leq n$, assume that the boundary $\partial \bD$ is divided into 
$2n$ arcs $a_1,\ldots, a_{k}, b_1, \ldots, b_{k}, a_{k+1}, \ldots a_{n},a_{k+1}^{-1}, \ldots a_{n}^{-1}$ 
such that the topological quotient  $\mathbb{P}^n_k := \bar\D\,/\!\sim\,$ with the 
equivalence relation \eqref{Psim} 
is homeomorphic to a closed non-orientable surface of Euler genus $n$. 
Then the C*-algebra 
\[ \label{Pq} 
\CPq := \{ f\hsp\in\hsp\T: \s(f)(x)\hsp= \hsp\s(f)(y),  \  \,\forall x,y\in\partial \bD \,\text{ such that } x\sim y \}  
\]
defines a  non-orientable closed quantum surface of Euler genus $n$. 
\end{defn} 
  
 Formally, we have defined a collection of quantum surfaces of the same genus. That is, all different arrangements 
 with the same number of oriented arcs that give classically the same surface define different quantum versions. 
 For instance, the two different 
 orders $a_1a_2a_1^{-1}a_2^{-1}\ldots a_{2g-1}a_{2g}a_{2g-1}^{-1}a_{2g}^{-1}$ and 
 $a_1a_2 \ldots  a_{2g-1}a_{2g}a_1^{-1}a_2^{-1}\ldots a_{2g-1}^{-1}a_{2g}^{-1}$ yield different subalgebras of $\T$, 
 but an orientable quantum surfaces of the same genus. 
 The classical cut-and-glue procedure for the classification of closed surfaces 
 does not apply here because the simple C*-algebra $\cC_0(\D_q):=\cK$ has no closed ideals, so it cannot be divided 
 into two pieces with a common boundary. In particular, we don't have any topological technique at our disposal 
 to prove that $\cC(\mathbb{P}^n_{k,q})$ and $\cC(\mathbb{P}^n_{k^\prime,q})$ are isomorphic for $k\neq k^\prime$. 
In fact, we shall show in Section \ref{Siso} that these C*-algebras are isomorphic  if and only if $k = k^\prime$. 
 On the other hand, all C*-algebras associated to the same {\it orientable} quantum surface are actually isomorphic.

 As illustrative examples and for the convenience of the reader, we will give an explicit description of a 
 closed quantum surfaces for each genus. For $g\in \N$, define $4g$ arcs on the circle $\bS$ by 
\begin{align*}
  a_k,\, a_k^{-1} : [0,1] \lra \bS,\ \  a_k(t) \hsp:=\hsp \e^{\pi\im\frac{k-1+t}{2g} }, \ \ 
  a_k^{-1}(t) \hsp:= \hsp\e^{\pi\im\frac{2g+ k-t}{2g} },\ \ k\hsp=\hsp 1,\ldots,2g. 
\end{align*}
Apparently, this arrangement differs from the usual ``normal form'' \cite{F}. Nevertheless 
$\bT^g:= \bar\D\,/\!\!\sim$ with the equivalence relation given in \eqref{Tsim} 
yields a closed orientable surface of genus $g$
and therefore \eqref{Tq} defines an orientable closed quantum surface of genus $g$. 
For an example of  a non-orientable closed quantum surface of Euler genus $n$, we may consider 
the $2n$ arcs 
\begin{align*}  
  a_k,\, b_k : \,[0,1]\, \lra\, \bS,\quad  a_k(t) := \e^{\pi\im\frac{k-1+t}{n} } , \quad b_k(t) := \e^{\pi\im\frac{-k+t}{n} },\quad k=1,\ldots,n. 
\end{align*}
In both cases, the sign before $t$ determines the orientation of the arcs. 
 
Note that,  as $\s$ is a *-homomorphism and therefore norm decreasing, 
our definitions yield indeed C*-subalgebras of $\T$. Moreover, 
$\cC_0(\D_q) \subset \CTq$ and $\cC_0(\D_q)\subset \CPq$ 
since $\cC_0(\D_q) := \cK = \ker\s$. On the boundary, the functions 
$\s(\CTq)\subset\CS$ and $\s(\CPq)\subset\CS$ 
do not separate the identified points along two equivalent arcs. For $\CTq$ and $\CPq$, 
it can be checked that all arcs start and end at the same point. 
Hence $\s(\CTq)$ and $\s(\CPq)$ separate the points of 
$2g$ and $n$ arcs, respectively, all starting and ending at the same point. 
As a consequence, 
$$
\s(\CTq) \cong C(\underset{k=1}{\overset{2g}{\vee}} \bS ), \quad 
\s(\CPq) \cong C(\underset{k=1}{\overset{n}{\vee}} \bS ), 
$$
where $\underset{k=1}{\overset{N}{\vee}} \bS$ denotes the wedge product of $N$ circles. 
In case of $\CSq$, the image of the symbol map yields only the continuous functions 
on a half circle $\bS_+:= \{x\in\bS : \mathrm{Im}(x) \geq 0\}$. 
This observation leads to the following C*-al\-ge\-bra extensions: 
\begin{align} \label{ext1}
 &\xymatrix{
 0\;\ar[r]&\; \cK\;\ar[r]^{\iota\ \ } &\; \CTq \;\ar[r]^{\s\ \ } &\;  \cC\big(\underset{k=1}{\overset{2g}{\vee}} \bS\big)\;\ar[r] &\ 0\, ,}& \\[-4pt] 
& \label{ext2} \xymatrix{
  0\;\ar[r]&\; \cK\;\ar[r]^{\iota\ \ } &\; \CPq \;\ar[r]^{\s\ \ } &\;  \cC\big(\underset{k=1}{\overset{n}{\vee}} \bS \big)\;\ar[r] &\ 0\, ,}&  \\[-4pt] 
  & \label{ext3}  \xymatrix{
 0\;\ar[r]&\; \cK\;\ar[r]^{\iota\ \ } &\ \CSq \;\ar[r]^{\s\ \ } &\  \cC( \bS_+)\ \ar[r] &\ 0\, .}&
\end{align}
The surjectivity can be verified by lifting a function $f\in\cC\big(\underset{k=1}{\overset{N}{\vee}} \bS \big)\subset \CS$ 
(or $f\in \cC(\bS_+)\subset \CS$) to a continuous function $\hat f\in \cC(\bD)$, $\hat f(r\e^{\im\theta}):=r\hs f(e^{\im\theta})$, 
and recalling that $\s(T_{\hat f})=f$. 

It is well known (see e.\,g.\ \cite{WO}) that a C*-extension 
gives rise to an isomorphic description 
as a pullback of C*-algebras via the Busby invariant. Let 
\[ \label{rohn} 
\rho_N \,:\, \bS\ \lra \ \big(\bS/\!\sim\!\big)\ \cong\ \underset{k=1}{\overset{N}{\vee}} \bS 
\]
denote the quotient map defined by restricting the quotients in \eqref{Tsim} and \eqref{Psim} to the boundary, 
where $N=2g$ and $N=n$, respectively. Then the inclusion $\cC\big(\underset{k=1}{\overset{N}{\vee}} \bS \big)\subset \CS$ 
corresponds to the pullback $\rho_N^* : \cC\big(\underset{k=1}{\overset{N}{\vee}} \bS \big)\ra \CS$ 
and the Busby invariant is determined by 
$$
\tau_N : \cC\big(\underset{k=1}{\overset{N}{\vee}} \bS \big) \lra \Cal=\cB/\cK, \quad \tau_N(f) = \s\big(T_{\widehat{\rho_N^*(f)}}\big),  
$$
where $\widehat{\rho_N^*(f)}$ stands for the extension of $\rho_N^*(f)\in \CS$ to the closed disk as described below \eqref{ext3}. 
By \cite[Prop.\ 3.2.11]{WO}, our closed quantum surfaces are naturally isomorphic to the pullback 
 \begin{equation}     \label{BI} 
\mbox{$\xymatrix@=2mm{& & \cB \underset{(\s,\tau_N)}{\oplus} \cC\big(\underset{k=1}{\overset{N}{\vee}} \bS \big) 
\ar[lld]^{\mathrm{pr}_1} \ar[rrd]_{\mathrm{pr}_2}  & &\\
\cB \ar@{>>}[drr]_{\s}& & & &\cC\big(\underset{k=1}{\overset{N}{\vee}} \bS \big) \,, \ar[dll]^{\tau_N}\\
&& \cB/\cK &&}$} 
\end{equation} 
where $\s:\cB\ra\cB/\cK$ denotes the quotient map. 
As the image of $\tau_N$ lies in $\s(\T) = \T/\cK \cong  \CS$, we obtain the same pullback C*-algebra by the reduced 
pullback diagram 
 \begin{equation} \label{CW} 
\mbox{$\xymatrix@=2mm{& & \cC(\bD_q) \underset{(\s,\rho_N^*)}{\oplus} \cC\big(\underset{k=1}{\overset{N}{\vee}} \bS \big) 
\ar[lld]^{\mathrm{pr}_1} \ar[rrd]_{\mathrm{pr}_2}  & &\\
\cC(\bD_q)  \ar@{>>}[drr]_{\s}& & & &\cC\big(\underset{k=1}{\overset{N}{\vee}} \bS \big) \,, \ar[dll]^{\rho_N^*}\\
&& \CS &&}$} 
\end{equation} 
where we made use of the quantum disk picture $\cC(\bD_q) :=\T$. 

The pullback diagram \eqref{CW} allows a nice interpretation of closed quantum surfaces 
as noncommutative CW-complexes \cite{qCW}. Classically, we may view \eqref{CW} as a dualization 
of the pushout diagrams 
\begin{equation} \label{CWD} 
\mbox{$\xymatrix@=3mm{& & \mathbb{T}^{g}  & &\\
\bD \ar[urr]& & & &\underset{k=1}{\overset{2g}{\vee}} \bS \,,\ar[ull]\\
&& \bS   \ar@{_{(}->}[ull]^{\iota}\ar[urr]_{\rho_{2g}} &&}$}
\qquad
\mbox{$\xymatrix@=3mm{& & \mathbb{P}^{n}  & &\\
\bD \ar[urr]& & & &\underset{k=1}{\overset{n}{\vee}} \bS \,,\ar[ull]\\
&& \bS  \ar@{_{(}->}[ull]^{\iota}\ar[urr]_{\rho_n} &&}$}
\end{equation}
where $\rho_{2g}$ in the left diagram and $\rho_{n}$ in the right diagram 
are given by the restriction to the boundary of the topological quotients 
defined by the equivalence relations in \eqref{Tsim} and \eqref{Psim}, respectively, 
and $\iota : \bS\cong \partial\bD \hookrightarrow \bD$ denotes the inclusion. 
Clearly, we may view $\underset{k=1}{\overset{N}{\vee}} \bS$ as a 1-skeleton obtained by attaching 
$N$ arcs to a 0-skeleton consisting of a single point. Then the diagrams in \eqref{CWD} amount 
to attaching a 2-cell to the 1-skeletons, so that \eqref{CW} becomes a dualized and  quantized version of it. 
A generalization of this construction to higher dimensions, including 
K-theoretic computations by using spectral sequences, will be given in \cite{SAW}.

 \section{Isomorphism classes of quantum surfaces}  \label{Siso} 
 
 In section we address the question of isomorphism classes of closed quantum surfaces. 
 As in Definition \ref{D1}, we will only allow the assignment of 
 $4g$ arcs for $\CTq$  and $2n$ arcs for $\CPq$ in such a way that 
 the construction yields the classical counterpart if the quantum disk gets replaced by the closed unit disk. 
Thus, in the orientable case, only arcs of opposite orientation are pairwise identified, and in the 
non-orientable case, there exists at least one pair of identified arcs having the same orientation 
on the boundary circle. 

For the purpose of applying Brown-Douglas-Fillmore theory \cite{BDF1,BDF2}, we will use the 
pullback diagram \eqref{BI} and characterize the C*-algebra extension \eqref{ext1} and \eqref{ext2} by a 
single, essentially normal generator.  To begin, we describe $\underset{k=1}{\overset{N}{\vee}} \bS$ 
homeomorphically as a compact subset in $\C$, for instance as a finite  Hawaiian earring:  
\[  \label{XN} 
\varphi_N: \underset{k=1}{\overset{N}{\vee}} \bS \,\overset{\cong}{\lra}\, 
X_N \,:=\, \underset{k=1}{\overset{N}{\cup}} \mathbb{S}^1_{\frac{k+1}{k}}(-\mfrac{1}{k}) 
= \underset{k=1}{\overset{N}{\cup}} \{ x\in \C : |x+\mfrac{1}{k}|= \mfrac{k+1}{k}\}. 
\] 
Here, $\mathbb{S}^1_{\frac{k+1}{k}}(-\mfrac{1}{k})$ stands for the circle with radius \mfrac{k+1}{k} and centre $-\mfrac{1}{k}\in\C$. 
All these circles have a common base point at $1\in\C$. Let $z:X_N\ra \C$, $z(x)=x$ denote 
the identity function. Clearly, $z$ separates the points of $X_N$, hence $z$, $\bar z$ and $1$ generate 
the C*-algebra $\cC(X_N)$ by the Stone--Weierstrass theorem. Thus the function $\zeta_N\in \cC(\bS)$, 
\[ \label{zN}
\zeta_N:= (\varphi_N \circ \rho_N)^*z : \bS \ra X_N\subset  \C
\] 
separates exactly the points of the arcs after the 
identification by $\rho_N$. Hence any function 
$h\in \cC\big(\underset{k=1}{\overset{N}{\vee}} \bS \big)\cong \cC(\bS/\!\!\sim) \subset \cC(\bS) \cong \s(\T)$ 
satisfying the same ``boundary conditions'' 
from Definition \ref{D1} as the function $\zeta_N$ can be 
approximated by polynomials in $\zeta_N$ and $\bar\zeta_N$. 
This means that any such $h$ can be approximated
by polynomials in $\s(T_{\widehat\zeta_N})$ and $\s(T_{\widehat\zeta_N}^*)$
with the usual extension 
of $\zeta_N\in\CS$ to a continuous function 
\[ \label{extend} 
\widehat\zeta_N\in\cC(\bD), \ \ 
\widehat\zeta_N(r\e^{\im \theta})=r\hs\zeta_N(e^{\im \theta}), \ \ r\in[0,1], \ \ \theta\in\R.
\]
Thus $\s(\CTq) \cong \cC( X_{2g})$ 
and $\s(\CPq) \cong \cC(X_n)$,  
so the closure of the *-al\-ge\-bra generated by 
$T_{\widehat\zeta_N}$, $T_{\widehat\zeta_N}^*$, $1$ and $\cK$ defines a C*-al\-ge\-bra extension which has a Busby invariant that 
is isomorphic to $\tau_N$ in \eqref{BI} via the homeomorphisms  
$X_N \cong \underset{k=1}{\overset{N}{\vee}} \bS \cong \big(\bS/\!\sim\!\big)$. 
Therefore the C*-al\-ge\-bra extensions are isomorphic, 
see Equation \eqref{Piso} below. By definition, the generator $T_{\widehat\zeta_N}$ yields 
an essentially normal operator with essential spectrum $X_N\subset \C$.  

Assume now that there exists an isomorphism of closed quantum surfaces 
$\a : \cC(\mathbb{M}_q) \ra \cC(\mathbb{M}_q^\prime)$, where 
$\mathbb{M}_q, \mathbb{M}_q^\prime \in\{\mathbb{T}^g_{q} , \mathbb{P}^n_{k,q}: g,n,k\in\N, \ k\leq n\}$. 
As all considered C*-algebras are subalgebras of $\T$ and contain 
the Jacobson radical $\cK=\ker(\s)$, we get 
an isomorphism $\a : \cK \ra \cK$. Any such isomorphism can be implemented by a 
unitary operator $U_\a\in\cB$ \cite[Remark 2.5.3]{HR}. Moreover, each isomorphism 
$\a : \cK \ra \cK$ has a unique extension to its multiplier algebra $\cB=\cM(\cK)$. 
Therefore the isomorphism $\a : \cC(\mathbb{M}_q) \ra \cC(\mathbb{M}_q^\prime)$ 
can be given by $\a(t) =U_\a t U_\a^*$ with a unique unitary operator $U_\a$. 
As a consequence, the C*-al\-ge\-bras 
$\cC(\mathbb{M}_q), \cC(\mathbb{M}_q^\prime) \subset \T$ are isomorphic 
if and only if there exist essentially normal operators 
$T\in\cC(\mathbb{M}_q)$ and $T^\prime \in \cC(\mathbb{M}_q^\prime)$ 
that generate together with $1$ and $\cK$ the corresponding C*-algebras 
and are unitarily equivalent up to a compact perturbation. 
The question of the existence of such a unitary equivalence 
is exactly the starting point of Brown-Douglas-Fillmore theory. 
This theory provides the principal tools for the classification of isomorphism classes 
in the next theorem.

\begin{thm}  \label{Tiso}
For $g\in \N$, let $\cC( \mathbb{T}^g_{q})$ denote an orientable closed quantum surface 
as defined in Definition \ref{D1}. Then $\cC( \mathbb{T}^g_{q})$ 
is isomorphic to the C*-al\-ge\-bra generated 
by $1$, $T_{g}$,  $T_{g}^*$ and the compact operators $\cK:=\cK(H_{g})$, where 
$$
H_{g} := \mplus{j=1}{2g} \lZ, \qquad 
T_{g}:=  \mplus{j=1}{2g} (\mfrac{j+1}{j}U\hsp - \hsp\mfrac{1}{j}), 
$$
and $U$ stands for the unitary bilateral shift on $\lZ$. 

Given $n,k\in\N$ with $k\leq n$, let $\cC( \mathbb{P}^n_{k,q})$ denote a 
non-orientable closed quantum surface from Definition \ref{D1}. 
Then $\cC( \mathbb{P}^n_{k,q})$ is isomorphic to the C*-al\-ge\-bra generated 
by $1$, $T_{n,k}$,  $T_{n,k}^*$ and the compact operators $\cK:=\cK(H_{n,k})$, where 
$$
H_{n,k} := \mplus{j=1}{k} \lN \oplus \mplus{j=k+1}{n} \lZ, \quad 
T_{n,k}:= \mplus{j=1}{k} (\mfrac{j+1}{j}S^2\hsp-\hsp\mfrac{1}{j}) \oplus\mplus{j=k+1}{n} (\mfrac{j+1}{j}U\hsp - \hsp\mfrac{1}{j}), 
$$
and $S$ stands for the unilateral shift from \eqref{S}. 

In particular, $\cC( \mathbb{T}^g_{q})$ and $\cC( \mathbb{T}^{g^\prime}_q)$ are isomorphic if and only if $g=g^\prime$. 
Furthermore, $\cC(\mathbb{P}^{n}_{k,q})$ and $\cC(\mathbb{P}^{n^\prime}_{k^\prime,q})$ are isomorphic 
 if and only if $n=n^\prime$ and $k=k^\prime$. Moreover, $\cC( \mathbb{T}^g_{q})$ is never isomorphic to $\cC(\mathbb{P}^{n}_{k,q})$, 
 and $\CSq$ is neither isomorphic to $\cC( \mathbb{T}^g_{q})$ nor to  $\cC(\mathbb{P}^{n}_{k,q})$.

\end{thm}    
\begin{proof}
First observe that, if we replace in the diagram \eqref{CW}  the C*-al\-ge\-bra 
$\cC\big(\underset{j=1}{\overset{N}{\vee}} \bS \big)$ by its isomorphic image 
$(\varphi^*_N)^{-1}: \cC\big(\underset{j=1}{\overset{N}{\vee}} \bS \big) \overset{\cong}{\lra}  \cC(X_N)$,  
with $\varphi_N$ given in \eqref{XN}, 
then we obtain isomorphic pullbacks 
\[ \label{Piso}
\phi_N : \cC(\bD_q)\! \underset{(\s,\rho_N^*\circ\varphi^*_N)}{\oplus}\!  \cC(X_N)
\overset{\cong}{\lra} \cC(\bD_q)\! \underset{(\s,\rho_N^* )}{\oplus}\! \cC\big(\underset{j=1}{\overset{N}{\vee}} \bS \big), 
\ \ \phi_N((t,f)) \! := \! (t, \varphi^*_N(f)). 
\]
Hence we can use the homeomorphism $\varphi_N$ from \eqref{XN} to order the circles from the wedge product 
in a certain ``normal form''. Suppose therefore without loss of generality that the closed quantum surface is defined
by the pairwise identification of $2N$ arcs such that the first $k$ circles correspond to pairs of arcs that had 
the same \emph{positive} orientation, and the remaining $N-k$ circles correspond to pairs of arcs that have been identified 
with opposite orientation. Here the homeomorphism $\varphi_N$ may be used to flip the orientation and to change the order of the arcs. 
As our assignment of arcs leads in the classical case to a closed surface, we get  
$\bar\D/\!\!\sim\ \cong \mathbb{P}^N_k$ similar to \eqref{TP} 
and $\mathbb{P}^N_0 = \mathbb{T}^{N/2}$ for $k=0$ by definition. 

The operators $\s(U), \s(S^2)\hsp \in\hsp \Cal\hsp=\hsp \cB/\cK$ have both the spectrum $\bS\hsp\subset\hsp \C$, 
therefore $\spec(T_{N,k}) = X_N \cong  \underset{j=1}{\overset{N}{\vee}} \bS$. 
However,  $\mathrm{Ind}(U)=0$ and $\mathrm{Ind}(S^2)=-2$. 
Further, $S^*\oplus S  + K \cong U$, where $K\in \cK(\lN\oplus\lN)\cong\cK(\lZ)$ maps $\mathrm{ker}(S^*)$ 
unitarily onto  $\mathrm{Im}(S)^\bot  \cong \mathrm{coker}(S)$. 
On the other hand,  it was explained in the beginning of this section that 
the C*-algebra of the closed quantum surface is generated by 
$T_{\widehat\zeta_N}$, $T_{\widehat\zeta_N}^*$, $1$ and $\cK$, where $\zeta_N$ 
has been defined in \eqref{zN} and $\widehat\zeta_N$ denotes its extension to the closed disk as in \eqref{extend}. 
So the proof of the theorem boils down to the question of 
when $T_{\widehat\zeta_N}$ is unitarily equivalent to a compact perturbation of $T_{N,k}$. 

Recall, e.\,g.\ from \cite[Section 16.2]{B}, that the essentially normal operators 
with essential spectrum $X_N\subset \C$ are classified, up to compact perturbations, by 
$K^1(X_N)\cong K_1(\cC(X_N))$.   
By \eqref{XN}, 
$K_1(\cC(X_N))\cong K_1( \cC(\underset{j=1}{\overset{N}{\vee}} \bS))$, and since the 
wedge sum of circles $\underset{j=1}{\overset{N}{\vee}} \bS$ can be obtained 
by the one-point compactification of $N$ open, disjoint intervals, we have that 
$K_1(\cC(X_N)) = \underset{j=1}{\overset{N}{\oplus}} \Z[u_j]$, 
where $u_j\in \cC(X_N)$ is any invertible function with winding number 1 (or -1) on the $j$-th circle 
and winding number 0 on all the others. 
Moreover, the winding number of an invertible  {Im}
function $\Phi\in\CS$ is related to the Fredholm index of the Toeplitz 
operator $T_{\widehat\Phi}$ by \eqref{wind}, and Brown-Douglas-Fillmore theory tells us that the Fredholm index 
is a principal obstruction for unitary equivalence of essentially normal operators up to compacts. 

Let $\zeta_N\in \CS$ and $\widehat \zeta_N \in \cC(\bD)$ be given by \eqref{zN} and \eqref{extend}, respectively. 
Then the essentially normal generator 
$T_{\widehat\zeta_N}\in\T$ of the corresponding closed quantum surface has essential spectrum 
$\mathrm{Im}(\s(T_{\widehat\zeta_N})) = \mathrm{Im}(\zeta_N) = X_N$.  
After applying $\varphi_N$ from \eqref{XN} to bring the circles into normal form, 
the function $\zeta_N$ from \eqref{zN} winds along an arc $a_j\subset \bS$ once around the 
circle $\mathbb{S}^1_{\frac{j+1}{j}}(-\mfrac{1}{j})$ in positive direction, and 
along the arc $a_j^{-1}\subset \bS$ once around the same 
circle $\mathbb{S}^1_{\frac{j+1}{j}}(-\mfrac{1}{j})$, but in negative direction. 
As the circles are ordered in normal form, and there are $2k$ arcs corresponding in the classical case to the 
connected sum $\mathbb{P}^1 \#\cdots  \#\mathbb{P}^1$ of $k$ projective spaces, the function $\zeta_N$ 
winds  exactly twice in positive direction around each of the first $k$ circles.  
For
the remaining $N- k$ circles, illustrated in \eqref{TP} as 
connected sum $\mathbb{T}^1 \#\cdots  \#\mathbb{T}^1$ of $(N - k)/2$ tori, 
any arc occurs also in the opposite direction, so the function  $\zeta_N$ has winding number 0 around each of 
these circles. From the classification of 
essentially normal operators by winding numbers in \cite[Theorem 16.2.1 and Example 16.2.4]{B}, 
together with the relation between winding numbers and the Fredholm index of shift operators in \eqref{Sind} and \eqref{wind}, 
it follows that the generator $T_{\widehat\zeta_N}$ is unitarily equivalent to a compact perturbation of 
$T_{N,k}$ defined in the theorem. Consequently the C*-algebra generated by 
$T_{\widehat \zeta_N}$, $T_{\widehat\zeta_N}^*$, $1$ and $\cK$ is isomorphic to the 
C*-algebra generated by $T_{N,k}$, $T_{N,k}^*$,  $1$ and $\cK$. 

Finally, two operators $T_{N,k}$ and $T_{N^\prime,k^\prime}$ 
are unitarily equivalent up to a perturbation by a compact operator if and only if they have the same essential 
spectrum, and $\s(T_{N,k})$ and $\s(T_{N^\prime,k^\prime})$ have the same winding numbers, i.e., $N=N^\prime$ and $k=k^\prime$. 
This implies the last claims of the theorem. 
\end{proof} 

Theorem \ref{Tiso} has two interesting consequences. 
First, we did not use in the proof the condition that the classical counterpart yields a closed surface. 
So there are assignments of arcs, always with starting and endpoint identified, that do not 
give rise to a 2-dimensional manifold in the classical case, but define a C*-algebra isomorphic 
to a closed quantum surface. Thus, on the one hand the Toeplitz quantization decreases degeneracy 
by distinguishing between $\cC(\mathbb{P}^{n}_{k,q})$ and $\cC(\mathbb{P}^{n}_{k^\prime,q})$ for $k\neq k^\prime$, 
and on the other hand it increases degeneracy by allowing for ``non-admissible'' 
prescriptions of arcs that do not even yield topological manifolds in the classical case. 

Second, there is an abuse of notation in Definition \ref{D1}. Equations \eqref{Tq} and  \eqref{Pq} define actually families 
of different C*-subalgebras  of $\T$, i.e., different arrangements yield different subalgebras.  
However, Theorem \ref{Tiso} shows that each admissible arrangement leads to a C*-algebra that is isomorphic to exactly one 
from Definition \ref{D1}.

\section{K-theory of closed quantum surfaces}  \label{SKT} 

In Section \ref{Sdef}, closed quantum surfaces were defined by analogy to the classical case. 
In this section, we will show that the topological invariants in the disguise of K-groups are not changed 
by the quantization process. A motivation for this fact was already given at the end of Section \ref{STA}. 

\begin{thm}  \label{TKG}
Let $\CTq$, $\CSq$ and $\CPq$ be as defined in Definition \ref{D1}. Then 
\begin{align*}
K_0(\CTq) &\cong \Z \oplus \Z,  & K_1(\CTq) &\cong\mplus{j=1}{2g} \Z, \\
K_0(\CSq) &\cong\Z \oplus \Z,  & K_1(\CSq) &\cong0, \\
K_0(\CPq) &\cong\Z_2 \oplus \Z,  & K_1(\CPq) &\cong\mplus{j=1}{n-1} \Z. 
\end{align*} 
In particular, all closed quantum surfaces from Definition \ref{D1} have the same K-groups 
as their classical counterparts. 
\end{thm} 

\begin{proof} The K-groups can easily be computed by applying the 6-term exact sequence of 
to the C*-algebra extensions  \eqref{ext1}--\eqref{ext3}: 
\begin{equation}  \label{stM}
\xymatrixcolsep{2pc}
\xymatrix{  
 K_0(\cK) \ \ar[r]^-{\iota_\ast} & 
  \ K_0 (\cC(\mathbb{M}_q))\   \ar[r]^-{{\s}_\ast} &
  \ K_0 (\cC(\underset{k=1}{\overset{N}{\vee}} \bS)) \ \ar[d]^{\mathrm{exp}}\\
 K_1 (\cC(\underset{k=1}{\overset{N}{\vee}} \bS))\ \ar[u]^{\ind} &
  \ K_1 (\cC(\mathbb{M}_q))\  \ar[l]_-{\ \ {\s}_\ast} &
\ K_1 (\cK) \ar[l]_-{\ \ \iota_\ast}  \,, 
  }
\end{equation} 
where $\mathbb{M}_q \in\{\cC(\mathbb{T}^g_{q}), \cC(\mathbb{P}^n_{k,q}), \cC(\mathbb{S}^2_q): g,n,k\in\N, k\leq n \}$. 
As discussed in Section \ref{STA}, $K_1 (\cK) = 0$ and $K_0 (\cK) = \Z[1\hsp-\hsp SS^*]$. 
Moreover, 
\begin{align*}
&K_0 (\cC(\underset{k=1}{\overset{N}{\vee}} \bS)) = K_0 \big((\underset{k=1}{\overset{N}{\oplus}}\cC_0(0,1))\dotplus \C\hs 1\big)
=  \big(\underset{k=1}{\overset{N}{\oplus}}\! K_0( \Sigma\C) \big) \oplus \Z [1] = \Z [1], \\ 
&K_1 (\cC(\underset{k=1}{\overset{N}{\vee}} \bS)) = K_1 \big((\underset{k=1}{\overset{N}{\oplus}}\cC_0(0,1))\dotplus \C\hs 1\big)
=  \underset{k=1}{\overset{N}{\oplus}}\! K_1(\cC_0(0,1))) = \mplus{j=1}{N} \Z, 
\end{align*} 
where $\cA\dotplus \C\hs 1$ means adjoining a unity to the non-unital C*-algebra $\cA$ 
and $\Sigma \cA$ denotes the suspension of $\cA$. Inserting these K-groups into \eqref{stM} yields 
\begin{equation}  \label{KM}
\xymatrixcolsep{2pc}
\xymatrix{  
\Z\ \ar[r]^-{\iota_\ast} & 
  \ K_0 (\cC(\mathbb{M}_q))\   \ar[r]^-{{\s}_\ast} &
  \ \Z[1] \ \ar[d]^{\mathrm{exp}}\\
 \Z^N\ \ar[u]^{\ind} &
  \ K_1 (\cC(\mathbb{M}_q))\  \ar@{_{(}->}[l]_-{\ \ {\s}_\ast} &
\ 0 \ar[l]_-{\ \ \iota_\ast}  \,, 
  }
\end{equation} 
Now $0 \ra \ker(\s_*) \ra K_0 (\cC(\mathbb{M}_q)) \overset{\s_*\ }{\lra} \Z [1]\ra 0$ is split exact 
with a splitting homomorphism given by $[1] \mapsto [1]$. Thus it follows from 
the exactness of \eqref{KM} that 
\[ \label{KM01} 
K_0(\cC(\mathbb{M}_q)) \cong \Z/\mathrm{Im}(\ind) \,\oplus \,\Z[1], \qquad K_1(\cC(\mathbb{M}_q)) \cong \mathrm{Ker}(\ind). 
\] 
Hence it remains to determine the index map $\ind: K_1 (\cC(\underset{k=1}{\overset{N}{\vee}} \bS)) \ra K_0(\cK)$. 

Recall that $\cC(\mathbb{M}_q) \subset \T$ and $\cC(\underset{j=1}{\overset{N}{\vee}} \bS) \subset \CS$ by Definition \ref{D1} 
and Equation~\eqref{Cex1}.  As explained at the end of Section \ref{STA},  
describing the index map amounts to lifting a unitary (matrix) $\Phi$  in 
$\cC(\underset{j=1}{\overset{N}{\vee}} \bS)$ to a Fredholm operator $F_\Phi$ in $\cC(\mathbb{M}_q) \subset \T$ and computing 
its Fredholm index $\mathrm{Ind}(F_\Phi)\in \Z \cong K_0(\cK)$. 
Moreover, the Fredholm index $\mathrm{Ind}(F_\Phi)$ coincides with the negative winding number  $- \wind[\Phi]$, see \eqref{wind}. 
We mentioned in the proof of Theorem \ref{Tiso} that 
$K_1(\cC(X_N)) = \underset{j=1}{\overset{N}{\oplus}} \Z[u_j]$, where $X_N\cong \underset{j=1}{\overset{N}{\vee}} \bS$ 
and $u_j$ is any invertible function that has winding number 1 (or -1) on the $j$-th circle 
and winding number 0 on all the others. Moreover, it was stated below \eqref{rohn} that 
the inclusion $\cC\big(\underset{j=1}{\overset{N}{\vee}} \bS \big)\subset \CS$ 
corresponds to the pullback $\rho_N^* : \cC\big(\underset{j=1}{\overset{N}{\vee}} \bS \big)\ra \CS$ 
with $\rho_N$ from \eqref{rohn}. 

Now let $\cC(\mathbb{M}_q) = \CTq$ and $N=2g$. Then, by \eqref{Tsim}, each arc, say $a_j$, occurs exactly once more 
with its negative orientation $a_j^{-1}$. As a consequence, the winding numbers of all invertible functions 
$\Phi \in  \cC\big(\underset{k=1}{\overset{2g}{\vee}} \bS \big)\subset \CS$ are 0, 
so $\ind \equiv 0$ and thus $K_0(\CTq) \cong \Z \oplus \,\Z[1]$ and $K_1(\cC(\mathbb{M}_q)) \cong \Z^{2g}$ 
by \eqref{KM} and \eqref{KM01}.   

Next we consider $\CPq$, $N=n$ and $1\leq k\leq n$. 
Assume that the circles of $\underset{j=1}{\overset{n}{\vee}} \bS$ are ordered in such a way that 
the first $k$ circles correspond to arcs that occur twice with the same orientation (i.e.\ $a_j(t) \sim b_j(t)$) 
and the other pairs with opposite orientations (i.e.\ $a_j(t) \sim a_j^{-1}(t)$). 
Therefore, an invertible function $u\in  \cC\big(\underset{k=1}{\overset{n}{\vee}} \bS \big)\subset \CS$ 
with winding number $m\in\Z$ along $a_j$ has also winding number $m$ along $b_j$ if $j\leq k$. 
On the other hand, if $j>k$, then a function with winding number $m\in\Z$ along $a_j$ 
will have winding number $-m$ along $a_j^{-1}$ so that these winding numbers add up to 0. 
For the generators $[u_j]$ of 
$K_1(\cC(\underset{j=1}{\overset{n}{\vee}} \bS)) \cong \underset{j=1}{\overset{N}{\oplus}} \Z[u_j]$
described above, we get 
$$
\ind[u_j] = 2, \ \ j\leq k, \qquad \ind[u_j] = 0, \ \ j>k, 
$$
so that $\ind(m_1,\ldots, m_n) = 2(m_1+\cdots + m_k)$  in the exact sequence \eqref{KM}. 
In particular, $\mathrm{Im}(\ind) = 2\Z$ and $\mathrm{Ker}(\ind) \cong \Z^{n-1}$ from which 
the result follows by \eqref{KM01}. 

In the case of $\CSq$ from \eqref{Sq}, a complex number $\e^{\pi\im t} \in \bS$ is 
identified with its complex conjugate $\e^{-\pi\im t} \in \bS$. The resulting quotient space $\bS/\!\sim$ 
is homeomorphic to an closed interval and therefore contractable. 
Replacing $\cC(\underset{j=1}{\overset{N}{\vee}} \bS)$ by $\C$ in \eqref{stM}, the lower row becomes 0 
and the upper row becomes exact, which yields the stated K-groups for $\CSq$. 

The last claim follows by comparing with the classical K-groups. 
\end{proof} 

For concrete calculations, it is convenient to have a suitable description of the generators of the K-groups. 
Let $u\in \CS$ be the identity function $u(\e^{\im t\theta}):= \e^{\im t\theta}$. Then the identity function 
$z\in \CD$, $z(r\e^{\im t\theta}):= r\e^{\im t\theta}$ is an extension of $u$ to the closed disk. 
Moreover, $[u]$ generates $K_1(\CS) \cong \Z[u]$. Set 
$$
P_{\text{Bott}}  := \begin{pmatrix} T_z T_z^*  &\!\!\! T_z \sqrt{1\!-\!\T_z^*T_z} \\[2pt]
                                                          \sqrt{1\!-\!T_z^*T_z}\,T_z^*   &\!\!\! 1\!-\! T_z^*T_z\end{pmatrix}
 =               \begin{pmatrix} T_z \\[2pt]
                                                          \sqrt{1\!-\!T_z^*T_z}\end{pmatrix}    
 \hsp\circ\hsp  \big(T_z^*,\sqrt{1\!-\!T_z^*T_z}\big)                                                                                               
 $$ 
Since the index map in \eqref{stX} is an isomorphism, it follows from \eqref{PTz} with $z$ instead of $\zeta$ 
that $\ind[u] = [P_{\text{Bott}}] -[1]$ 
generates $K_0(\cK)$. Note that $ [P_{\text{Bott}}] -[1]$ is never in the image of the index map from \eqref{stM}. 
Hence, for any closed quantum surface from Definition \ref{D1}, the $K_0$-group is generated 
by $[1]$ and $ [P_{\text{Bott}}]$. However, in the case of non-orientable quantum surfaces, 
we have the relation $2( [P_{\text{Bott}}] -[1]) =0$ as this element belongs to the image of the index map. 
On the other hand, we can lift $u\in \CS$ to the shift operator $S\in \T$, see \eqref{sS}. 
As described in the paragraph before \eqref{Sind}, $\ind[u] = -[1\!-\!SS^*]$, thus the relation 
$ [1] -[P_{\text{Bott}}] = [1\!-\!SS^*]$ holds in $K_0(\cK)$.

To describe the generators of the $K_1$-groups, consider (non-unitary) generators $[v_j]$ of 
$K_1 (\cC(\underset{k=1}{\overset{N}{\vee}} \bS)) \cong  K_1 (\cC(X_{N}))$ with winding number 1 
along the $j$-th circle and winding number 0 along the others.  More explicitly, set 
\[   \label{uj}
 v_{j}(x)= x, \ \, x \in \mathbb{S}^1_{\frac{j+1}{j}}(-\mfrac{1}{j}) \subset X_{N} \cong \underset{k=1}{\overset{N}{\vee}} \bS, \ \ 
 v_{j}(x)= 1\ \text{otherwise}, 
\]
and $u_j:=\rho_N^*(v_j) \in \CS$ with $\rho_N$ from \eqref{rohn}. 
Let $\widehat u_j\in \CD$ denote the extension of $u_j$ to the closed disk as given in \eqref{extend}. 
In the proof of the last theorem,  
we have seen that $\ind(T_{\widehat{u}_j}) = - \wind(u_j)=0$ in the orientable case, thus 
$\dim(\mathrm{Ker}(T_{\widehat{u}_j})) = \dim(\mathrm{Coker}(T_{\widehat{u}_j}))$. 
Choosing a compact isometry $K_j$ between $\mathrm{Ker}(T_{\widehat{u}_j})$ and $\mathrm{Im}(T_{\widehat{u}_j})^\bot$ 
and defining $T_j:= T_{\widehat{u}_j} +K_j$, we get an invertible operator in $\CTq$ 
such that $\s(T_j) = u_j$. Hence  $T_j$, or equivalently $U_j:= T_j\hs |T_j|^{-1}$, $j=1,\ldots,2g$, generate 
$K_1(\CTq)\cong \s_*(K_1(\CTq))\cong K_1 (\cC(\underset{k=1}{\overset{2g}{\vee}} \bS)) \cong  K_1 (\cC(X_{2g}))$. 
Under the isomorphism from Theorem~\ref{Tiso},  
the operator $T_j$ corresponds to a compact perturbation of 
$$ 
 \id\oplus \cdots \oplus \id \oplus  (\mfrac{j+1}{j}U\hsp - \hsp\mfrac{1}{j})\oplus \id\oplus \cdots \oplus \id 
 \,\in\, \cB(H_g) 
$$
as both have the same essential spectrum and the same winding numbers. 
 
In the \hsp non-oriented case, we consider the functions $v_{1,j}$ on $X_n\cong \underset{k=1}{\overset{n}{\vee}} \bS$ 
given by 
$$
v_{1,j}(x)\hsp  =\hsp  \bar x, \ x\hsp \in\hsp  \mathbb{S}^1_{2}(-1), \ \, 
v_{1,j}(x)\hsp =\hsp  x, \ x \hsp \in\hsp  \mathbb{S}^1_{\frac{j+1}{j}}(-\mfrac{1}{j}), \ \, 
v_{1,j}(x)\hsp  =\hsp  1\ \text{otherwise}, 
$$
for $j=2,\ldots, k$,
and $v_{1,j}:=v_j$ for $j>k$ with $v_j$ from \eqref{uj}. 
As before, let $u_{1,j}:= \rho^*_n(v_{1,j})$ denote  its pullback to $\CS$ 
and write $\widehat{u_{1,j}}\in\CD$ for its extension to the closed disk as in \eqref{extend}.  
Note that $u_{1,j}\in\CS$ has winding number 0 for all $j=2,\ldots, n$ so that $\ind[v_{1,j}]=\ind[u_{1,j}]=0$. 
As a consequence, $\mathrm{Ind}(T_{\widehat{u_{1,j}}}) =0$, or equivalently,  
$\dim(\mathrm{Ker}(T_{\widehat{u_{1,j}}})) = \dim(\mathrm{Coker}(T_{\widehat{u_{1,j}}}))$. 
Choosing compact isometries  $G_j $ between $\mathrm{Ker}(T_{\widehat{u_{1,j}}})$ and $\mathrm{Im}(T_{\widehat{u_{1,j}}})^\bot$, 
the operators $R_i:= T_{\widehat{u_{1,i+1}}} + G_{i+1}$, $i=1, \ldots, n-1$, become invertible 
in $\CPq$ and $\s(R_i) = u_{1,i+1}$. Thus $R_i$,  or equivalently $V_i:= R_i\hs |R_i|^{-1}$,  
defines an  element in $K_1(\CPq)$. 
Comparing the function $ u_{1,i+1}$ with the generators $[v_j]$ of 
$K_1 (\cC(\underset{k=1}{\overset{n}{\vee}} \bS)) \cong  K_1 (\cC(X_{n}))$ 
given in \eqref{uj}, we see that $\s_*([R_i])= [u_{1,i+1}]  = - [v_1] + [v_{i+1}]$ 
by counting the winding numbers along circles. In particular, $\s_*([R_i]) \in \mathrm{Ker}(\ind)$. 
Moreover, $ [v_{i+1}] - [v_1] $, $i=1,\ldots,n-1$, generate $\mathrm{Ker}(\ind)\cong \Z^{n-1}$. 
Therefore $[R_1]=[V_1], \ldots, [R_{n-1}]=[V_{n-1}]$ generate $K_1(\CPq)\cong\mathrm{Ker}(\ind)$. 
Finally, under the isomorphism from Theorem~\ref{Tiso}, the operator $R_i$ is 
a compact perturbation of 
\begin{align*}
 &(2\hs S^{*2}\hsp -\hsp 1) \hsp\oplus\hsp  \id \hsp\oplus\hsp  \cdots \hsp\oplus\hsp  \id \hsp\oplus\hsp   
 (\mfrac{i+2}{i+1}S^2\hsp - \hsp\mfrac{1}{i+1})\hsp\oplus\hsp  \id\hsp\oplus\hsp  \cdots \hsp\oplus\hsp  \id  
 \hs\in\hs\cB(H_{n,k}), \, \ i< k, \\
 &\id\hsp\oplus\hsp  \cdots \hsp\oplus\hsp  \id \hsp\oplus\hsp   
 (\mfrac{i+2}{i+1}U\hsp - \hsp\mfrac{1}{i+1})\hsp\oplus\hsp  \id\hsp\oplus\hsp  \cdots \hsp\oplus\hsp  \id
 \in\cB(H_{n,k}), \ \,i\hsp \geq\hsp k,
\end{align*} 
as these operators (for $i$ fixed) have the same essential spectrum and the same winding numbers. 

  \section*{Acknowledgements} 
 The authors acknowledge support from CIC-UMSNH, 
 from the CONACYT project A1-S-46784 ``Grupos cuánticos y geometría no conmutativa". 
 This research is part of the EU Staff Exchange project 101086394 "Operator Algebras That One Can See". 
 It was partially supported by the University of Warsaw Thematic Research Programme "Quantum Symmetries". 



\end{document}